\documentclass[11pt,leqno]{article}
\usepackage{amsthm,amsfonts,amssymb,amsmath,oldgerm}
\usepackage{epsfig}
\numberwithin{equation}{section}
\usepackage[thinlines]{easybmat}

\def\eps{\varepsilon }

\def\eps{\varepsilon}

\newcommand\br{\begin{remark}}
\newcommand\er{\end{remark}}
\newcommand\bp{\begin{pmatrix}}
\newcommand\ep{\end{pmatrix}}
\newcommand\be{\begin{equation}}
\newcommand\ee{\end{equation}}
\newcommand\ba{\begin{equation}\begin{aligned}}
\newcommand\ea{\end{aligned}\end{equation}}

\newcommand{\bap}{\begin{app}}
\newcommand{\eap}{\end{app}}
\newcommand{\begs}{\begin{exams}}
\newcommand{\eegs}{\end{exams}}
\newcommand{\beg}{\begin{example}}
\newcommand{\eeg}{\end{exaplem}}
\newcommand{\bpr}{\begin{proposition}}
\newcommand{\epr}{\end{proposition}}
\newcommand{\bt}{\begin{theorem}}
\newcommand{\et}{\end{theorem}}
\newcommand{\bc}{\begin{corollary}}
\newcommand{\ec}{\end{corollary}}
\newcommand{\bl}{\begin{lemma}}
\newcommand{\el}{\end{lemma}}
\newcommand{\bd}{\begin{definition}}
\newcommand{\ed}{\end{definition}}
\newcommand{\brs}{\begin{remarks}}
\newcommand{\ers}{\end{remarks}}

\newtheorem{theo}{Theorem}[section]
\newtheorem{prop}[theo]{Proposition}
\newtheorem{cor}[theo]{Corollary}
\newtheorem{lem}[theo]{Lemma}

\newtheorem{exams}[theo]{Examples}

\numberwithin{equation}{section}

\newtheorem{theorem}{Theorem}[section]
\newtheorem{proposition}[theorem]{Proposition}
\newtheorem{corollary}[theorem]{Corollary}
\newtheorem{lemma}[theorem]{Lemma}
\newtheorem{definition}[theorem]{Definition}

\newtheorem{example}[theorem]{Example}
\newtheorem{remark}[theorem]{Remark}

\newcommand{\RM}{\mathbb{R}}
\newcommand{\ZM}{\mathbb{Z}}

\title{
Nonlinear stability of spatially-periodic traveling-wave solutions of
systems of reaction diffusion equations
}

\author{\sc \small
Mathew A. Johnson\thanks{Indiana University, Bloomington, IN 47405;
matjohn@indiana.edu: Research of M.J. was partially supported by an NSF Postdoctoral Fellowship under NSF grant DMS-0902192.}
~~~~~
Kevin Zumbrun\thanks{Indiana University, Bloomington, IN 47405;
kzumbrun@indiana.edu: Research of K.Z. was partially supported
under NSF grants no. DMS-0300487 and DMS-0801745.
 }}
\begin{document}

\maketitle


\begin{center}
\end{center}


\begin{abstract}
Using spatial domain techniques developed by the authors and Myunghyun
Oh in the context of parabolic conservation laws,
we establish under a natural set of spectral stability conditions
nonlinear asymptotic stability with decay at Gaussian rate of
spatially periodic traveling-waves of systems of reaction diffusion
equations.
In the case that wave-speed is identically zero for all periodic solutions,
we recover and slightly sharpen a well-known result of Schneider obtained by renormalization/Bloch transform techniques; by the same arguments, we are able to treat
the open case of nonzero wave-speeds to which Schneider's
renormalization techniques do not appear to apply.
\end{abstract}

\section{Introduction}

In this paper, we study the nonlinear stability
with respect to spatially localized perturbations
of spatially periodic
traveling-wave solutions $u(x,t)=\bar u(x-ct)$ of a
system of reaction diffusion equations of form
$u_t=u_{xx}+f(u)$, where $(x,t)\in\RM\times\RM^+$, $u\in\RM^n$,
and $f:\RM^n\to\RM^n$ is sufficiently smooth:
equivalently, spatially periodic standing-wave solutions
$u(x,t)=\bar u(x) $ of
\begin{equation}\label{eqn:rd}
u_t-cu_x=u_{xx}+f(u).
\end{equation}
For the Allen--Cahn (variational) case $f(u)=dF(u)$,
the traveling-wave ODE becomes
$-cu'=u''+dF(u)$, hence for $c\ne 0$ either increases or
decreases a Hamiltonian
$$
H:=\frac{|u'|^2}{2} + F(u),$$
from which it follows that periodics exist only for speed $c=0$.
In this zero-speed case (not necessarily originating from
variational form), Schneider \cite{S} proved nonlinear stability
with decay at Gaussian (diffusive) rate by a combination of
weighted energy estimates, renormalization
techniques, and a spectacular nonlinear cancellation estimate,
all carried out in the Bloch frequency domain.
In the process, he showed that behavior is purely diffusive, with no
associated convection.
See also \cite{U,MSU} and references therein.

However, as described in \cite{DSSS}, there exist many other cases for
which there exist spatially periodic solutions of {\it varying speed},
in which situation one expects asymptotic behavior driven by
nonzero convection as well as diffusion.
In this case, the renormalization argument of \cite{S} does not
directly apply.

Meanwhile, motivated strongly by the work of Schneider, we and co-authors have carried out by somewhat different methods stability of
spatially periodic solutions of
systems of parabolic conservation laws $u_t + h(u)_x=u_{xx}$, for which
convection plays a major role \cite{OZ,JZ}.
Natural questions are (i) whether these alternative methods might be used to reproduce the original results of Schneider in the zero-speed case, and
(ii) whether, more, they might be able to treat the nonzero-speed (convective)
case left open up to now.

In this paper, we answer both questions in the affirmative, establishing
stability and decay at Gaussian rate for the general case, with no condition
on the wave-speed; see Theorem \ref{main} in Section \ref{s:nstab}.
Moreover, our nonlinear iteration method,
based on spatial rather than (Bloch) frequency
domain, permits an extremely brief and simple proof yielding new insight
even in the zero-speed case treated previously by Schneider.
\medskip

{\bf Note:}
We have been informed by B\"jorn Sandstede that
a similar result has been obtained by different means in \cite{SSSU}
using a nonlinear decomposition of phase and amplitude variables as
in \cite{DSSS}, accomodating also nonlocalized perturbations
in the phase.\footnote{
In the notation of Thm. \ref{main},
data $\tilde u_0$ with
$\|\tilde u_0(x+\psi_0)-\bar u\|_{L^1(\RM)\cap H^K(\RM)}|$,
$|\partial_x \psi_0|_{H^K(\RM)}$
sufficiently small.
		}

\medskip

{\bf Acknowledgement.} Thanks to Blake Temple, for his encouragement and
interest in the problem.  Also, we are greatful to B\"jorn Sandstede for describing the setting of the nonzero speed problem and for pointing out the reference \cite{SSSU}.

\section{Existence and Spectral Stability Assumptions}

Any standing wave solution of \eqref{eqn:rd} clearly must be a solution of the ordinary differential equation
\begin{equation}
u_{xx}+cu_x-f(u)=0,\label{eqn:peq}
\end{equation}
which is commonly referred to as the profile equation or as the traveling wave ODE corresponding to the original system.
The existence of periodic orbits of \eqref{eqn:peq} is trivial in the case $c=0$ where the equation is clearly Hamiltonian and hence
can be directly integrated by quadrature.  When the wave-speed is non-zero, however, the existence of periodic orbits is more
delicate but still straightforward: it can be treated via ODE/implicit function theorem techniques as familiar in the conservation law
case.  Indeed, writing \eqref{eqn:peq} as a $2n \times 2n$ system with $2n$ constraints (periodicity)
and two extra parameters (speed and period) yields, generically, a two-dimensional solution set.  The techniques
presented in this paper apply regardless of whether $c$ is non-zero and hence we study \eqref{eqn:peq} for a
general wave-speed $c\in\RM$.

Throughout our analysis, we assume the existence of an $X$-periodic
solution $\bar{u}(x)$ of \eqref{eqn:peq}.  It follows that generically one expects the periodic orbits to form a two-parameter
family of solutions of the profile equation, parameterized by the wave-speed $c$ and a translation mode $x_0$.  More precisely, we
make the following generic assumptions:
\begin{enumerate}
  \item[(H1)] $f\in C^K(\RM)$ for some $K\geq 2$.
  \item[(H2)]  The set of periodic solutions of \eqref{eqn:peq} in the vicinity of $\bar{u}$ forms a smooth $2$-dimensional
   manifold $\{\bar{u}^{x_0,s}(x-st+x_0)\}$ with $x_0,s\in\RM$
\end{enumerate}

We begin our
stability analysis by considering
the linearization of \eqref{eqn:rd} about the fixed
periodic standing-wave solution $\bar u$.
Without loss of generality, we assume that $\bar u$ is $1$-periodic, i.e.
that $\bar u(x+1)=\bar u(x)$ for all $x\in\RM$.  Considering nearby solutions of the form
\[
\bar u(x)+\eps v(x,t)+\mathcal{O}(\eps^2),
\]
where $|\eps|\ll 1$ and $v(\cdot,t)\in L^2(\RM)$, corresponding to spatially localized perturbations,
we see that $v$ satisfies the linear equation
\[
v_t=v_{xx}+df(\bar u)v.
\]
Since this equation is autonomous in time, we may seek separated solutions of the form
\[
v(x,t)=e^{\mu t}v(x)
\]
which readily yields the spectral problem
\begin{equation}
\mu v=L[u]v:=\left(\partial_x^2+ c\partial_x + df(\bar u)\right)v\label{eqn:spec}
\end{equation}
considered on the real Hilbert space $L^2(\RM)$, where here we are considering the linear operator $L[u]$ as having dense
domain in $H^2(\RM)$.  As coefficients of $L[u]$ are $1$-periodic, Floquet theory implies the $L^2$ spectrum is purely continuous
and corresponds to the union of the $L^\infty$ eigenvalues corresponding to considering the linearized operator
with boundary conditions $v(x+T)=e^{i\kappa}v(x)$ for all $x\in\RM$, where $\kappa\in[-\pi,\pi]$ is referred to as the Floquet
exponent and is uniquely defined mod $2\pi$.  In particular, $\mu\in\sigma(L[u])$ if and only if the spatially periodic
spectral problem \eqref{eqn:spec} admits a bounded eigenfunction of the form
\begin{equation}
v(x)=e^{i\xi x}w(x)\label{eqn:ef}
\end{equation}
where $w(x+1)=w(x)$.

Substitution of the Ansatz \eqref{eqn:ef} into \eqref{eqn:spec} motivates the use of the Fourier-Bloch decomposition
of the spectral problem.  To this end, we follow \cite{G,S} and define the one-parameter family of linear operators,
referred to as Bloch operators, by
\[
L_{\xi}[u]:=e^{-i\xi x}L[u]e^{i\xi x},~~\xi\in[-\pi,\pi],
\]
operating on $L^2_{\rm per}([0,1])$.  The $L^2$ spectrum of the linearized operator $L[u]$ is readily seen to be given by
the union of the spectra of the Bloch operators.  By continuity of the spectrum on the Floquet parameter $\xi$, and the discreteness
of the spectrum of the elliptic operator $L[u]$ on the compact domain $[0,1]$, it follows that the spectra of $L[u]$ may be described
as  the union of countably many continuous surfaces $\mu(\xi)$.

Continuing with this functional setup, we recall that any function $g\in L^2(\RM)$ admits an inverse Bloch-Fourier representation
\[
g(x)=\frac{1}{2\pi}\int_{-\pi}^\pi e^{i\xi x}\hat{g}(\xi,x)d\xi
\]
where $\hat{g}(\xi,x)=\sum_{j\in\ZM}e^{2\pi ijx}\hat{g}(\xi+2\pi j)$ is a $1$-periodic functions of $x$,
and $\hat{g}(\cdot)$ denoting with slight abuse of notation the usual Fourier transform of the function $g$ in
the spatial variable $x$.  Indeed, using the Fourier transform we have
\begin{align*}
2\pi g(x)&=\int_{-\infty}^\infty e^{i\xi x}\hat{g}(\xi)d\xi
=\sum_{j\in\ZM}\int_{-\pi}^{\pi}e^{i(\xi+2\pi j)x}\hat{g}(\xi+2\pi j)d\xi
=\int_{-\pi}^\pi e^{i\xi x}\hat{g}(\xi,x)d\xi,
\end{align*}
where the summation and integral can be interchanged for Schwarz functions $g$.
By similar computations, it is also seen by Parseval's identity that the Bloch-Fourier
transform $g(x)\to\hat{g}(\xi,x)$ is an isometry
of $L^2(\RM)$, i.e.
\begin{equation}
\|g\|_{L^2(\RM)}=\int_{-\pi}^\pi\int_0^1|\hat{g}(\xi,x)|^2dx~d\xi=:\|\hat{g}\|_{L^2(\xi;L^2(x))}.\label{eqn:parseval}
\end{equation}
Moreover, the Bloch-Fourier transform diagonalizes the periodic-coefficient operator $L[u]$, yielding the
inverse Bloch-Fourier transform representation
\begin{equation}
e^{L[u]t}g(x)=\frac{1}{2\pi}\int_{-\pi}^\pi e^{i\xi x}e^{L_{\xi}[u]t}\hat{g}(\xi,x)d\xi\label{eqn:invbf}
\end{equation}
relating the behavior of the linearized system to that of the diagonal operators $L_{\xi}[u]$.

We now discuss the spectral stability of the underlying solution $u(x)$ in more detail.
To begin, notice that by the translation invariance of \eqref{eqn:rd} the function $u'(x)$ is a $1$-periodic solution of the
differential equation $L[u]v=0$.  Hence, it follows that $\mu=0$ is an eigenvalue of the Bloch operator $L_0[u]$.
From the structure of the equation, it is natural to assume that translation generates the only null direction of the
operator $L_0[u]$.  Indeed, notice by $(H2)$ and the secular dependence of the periodic on the wave-speed $c$, that variation in the
translation direction $x_0$ is the only generator of $1$-periodic null-directions of the linearized operator given by Noether's theorem.  With this in mind,
following \cite{S}, we make the following natural spectral stability assumptions:
\begin{enumerate}
  \item[(D1)] $\mu=0$ is a simple eigenvalue of $L_0[u]$.
(Recall that $\xi=0$ corresponds to co-periodic perturbations.)
  \item[(D2)]
$\Re \sigma(L_{\xi}[u])\leq-\theta|\xi|^2$ for some constant $\theta>0$.
\end{enumerate}
\noindent Assumptions (D1)-(D2) correspond to ``dissipativity" of the large-time bahavior of the
linearized system and are often referred to as \emph{strong}
or \emph{diffusive} spectral stability assumptions \cite{JZ,OZ,S}.

\br
\textup{
By standard spectral perturbation theory \cite{K}, (D1) implies
that the eigenvalue $\mu(\xi)$ bifurcating from $\mu=0$ at $\xi=0$
is analytic at $\xi=0$, with
$\mu(\xi)=\mu_1\xi+ \mu_2\xi^2 + \mathcal{O}(|\xi|^3)$,
from which we find from the necessary stability condition
$\Re \mu(\xi)\le 0$ that $\Re \mu_1=0$ and $\Re \mu_2\le 0$.
Assumption (D2) thus amounts to the nondegeneracy condition
$\Re \mu_2\ne 0$ together with the strict stability condition
$\Re \sigma L_\xi<0$ for $\xi\ne 0$.
Condition (D2) is never satisfied in the scalar case $n=1$, by
Sturm--Liouville considerations, hence the emphasis in the
title on {\rm systems} of reaction diffusion equations.
}
\er

The goal of our analysis is to prove that the above spectral stability assumptions
imply nonlinear $L^1\bigcap H^K\to H^K$ stability
of the underlying periodic traveling wave.  To this end, we use this spectral information to obtain bounds on the linearized solution operator
$e^{L[u]t}$.  As we will see, assumptions (D1)-(D2), along with (H1)-(H2),
imply that the solution operator decays
(in a suitable sense) polynomially in time at a fast enough rate to prove the nonlinear stability of the underlying spatially periodic solution.
These bounds are established in the next section, after which we present our nonlinear iteration scheme.

\section{Linear Estimates}

In this section, we make use of the spectral stability assumptions of the previous section in order to prove bounds on the solution
operator $S(t):=e^{L[u]t}$.  These linearized estimates form the crux of the nonlinear analysis, presented in the next section.
To begin, notice that by standard spectral perturbation theory
\cite{K}, assumption $(D1)$ implies that the total eigenprojection $P(\xi)$
onto the eigenspace of $L_{\xi}[u]$ associated with the eigenvalue $\mu(\xi)$ bifurcating from the $(\xi,\mu(\xi))=(0,0)$
state is well defined and analytic in $\xi$ for $\xi$ sufficiently small, since the discreteness of the spectrum
of $L_{\xi}[u]$ implies that the eigenvalue $\mu(\xi)$ is separated at $\xi=0$ from the remainder of the spectrum of $L_0[u]$.
In particular, there exists an eigenfunction $q(x,\xi)$ bifurcating from the $q(x,0)=u'(x)$ state defined for $|\xi|\ll 1$
such that
\begin{equation}
L_{\xi}[u]q(x,\xi)=\mu(\xi)q(x,\xi)\label{eqn:bif}
\end{equation}
where, by assumption $(D2)$, the function $\mu(\xi)$ satisfies the estimate
\begin{equation}
\Re(\mu(\xi))\leq-\theta|\xi|^2\label{eqn:d3}
\end{equation}
for some constant $\theta>0$.

Our strategy is to treat the high- and low-frequency parts of the full solution operator $S(t)$ separately since, as is typical,
the low-frequency analysis is considerably more delicate than the corresponding high-frequency analysis.
To this end, we introduce a smooth cut off function $\phi(\xi)$ such that
\[
\phi(\xi)=\left\{
            \begin{array}{ll}
              1, & \textrm{ if }|\xi|\leq \eps \\
              0, & \textrm{ if }|\xi|\geq 2\eps,
            \end{array}
          \right.
\]
where $\eps>0$ is a sufficiently small parameter, and we split the solution operator $S(t)$ into its low-frequency part
\begin{equation}
S^{I}(t)g(x):=\frac{1}{2\pi}\int_{-\pi}^\pi e^{i\xi x}\phi(\xi)P(\xi)e^{L_{\xi}[u]t}\hat{g}(\xi,x)d\xi\label{eqn:lfsolnop}
\end{equation}
and high frequency part
\[
S^{II}(t)g(x):=\frac{1}{2\pi}\int_{-\pi}^\pi e^{i\xi x}\left(1-\phi(\xi)P(\xi)\right)e^{L_{\xi}[u]t}\hat{g}(\xi,x)d\xi,
\]
by which one may readily check that $S(t)g=\left(S^{I}(t)+S^{II}(t)\right)g$ by \eqref{eqn:invbf}.
As the low-frequency analysis is more delicate, we begin by deriving $L^2\to L^p$ bounds on $S^{II}(t)$.


Using the fact that $L_{\xi}[u]$ is a sectorial operator, and the spectral separation of $\mu(\xi)$ from the
remaining spectrum of $L_{\xi}[u]$, standard semigroup theory\cite{He,Pa} trivially implies that the bounds
\begin{equation}\label{semigrpbds}
\begin{aligned}
\left\|e^{L_{\xi}[u]t}\left(1-\phi(\xi)P(\xi)\right)g\right\|_{L^2([0,1])}&\lesssim e^{-\theta t}\|g\|_{L^2([0,1])},\\
\left\|e^{L_{\xi}[u]t}\left(1-\phi(\xi)P(\xi)\right)\partial_x g\right\|_{L^2([0,1])}&\lesssim t^{-1/2} e^{-\theta t}\|g\|_{L^2([0,1])},\\
\left\|\partial_x e^{L_{\xi}[u]t}\left(1-\phi(\xi)P(\xi)\right)g\right\|_{L^2([0,1])}&\lesssim t^{-1/2} e^{-\theta t}\|g\|_{L^2([0,1])},\\
\end{aligned}
\end{equation}
for all $t>0$ and 
some constant $\theta>0$.  Together with \eqref{eqn:parseval},
this yields immediately the following estimate.

\begin{prop}\label{prop:hfbds}
Under assumptions (H1)-(H2) and (D1)-(D2), there exists a constant $\theta>0$ such that for all $2\leq p\leq\infty$
and $t>0$ we have
\begin{align*}
\|S^{II}(t)g\|_{L^p(\RM)}&\lesssim t^{-\frac{1}{2}(1/2-1/p)}e^{-\theta t}\|g\|_{L^2(\RM)}.
\end{align*}
\end{prop}

\begin{proof}
First, notice the bounds in \eqref{semigrpbds} and the triangle inequality imply that
\begin{align*}
\|\partial_x^m S^{II}(t)g\|_{L^2(\RM)}&\lesssim\int_{-\pi}^\pi \left\|\partial_x^m\left(1-\phi(\xi)P(\xi)\right)e^{L_{\xi}[u]t}\hat{g}(\xi,\cdot)\right\|_{L^2(x;[0,1])}d\xi\\
&\lesssim t^{-m/2} e^{-\theta t}\int_{-\pi}^\pi\|\hat{g}(\xi,\cdot)\|_{L^2(x;[0,1])}d\xi
\\ &= t^{-m/2}e^{-\theta t}\|g\|_{L^2(\RM)}
\end{align*}
for either $m=0$ or $m=1$, where the final equality is justified by \eqref{eqn:parseval}, thus justifying the first claim by taking $m=0$.

To prove the second inequality, note when $p=\infty$ Sobolev embedding and the above $L^2\to L^2$ bound implies that
\begin{align*}
\|S^{II}(t)g\|_{L^\infty(\RM)}&\lesssim\left(\| S^{II}(t)g\|_{L^2(\RM)}\cdot\|\partial_x S^{II}(t)g\|_{L^2(x;\RM)}\right)^{1/2}
\lesssim t^{-1/4} e^{-\theta t}\|g\|_{L^2(\RM)}.
\end{align*}
The result for general $2\leq p\leq \infty$ now follows by $L^p$ interpolation.
\end{proof}


In order to analyze the low-frequency part of the solution operator $S(t)$, we find it convenient to introduce the Green kernel
\[
G^I(x,t;y):=S^{I}(t)\delta_y(x)
\]
associated with $S^I$, and
\[
\left[G^I_\xi(x,t;y)\right]:=\phi(\xi)P(\xi)e^{L_{\xi}[u]t}\left[\delta_y(x)\right]
\]
the corresponding kernel appearing within the Bloch-Fourier representation of $G^I$, where the brackets $[\cdot]$ denotes the periodic
extensions of the given function onto the whole real line.  Our first goal is to provide a useful representation for $G^I$ which incorporates
the spectral properties (D1)-(D2) of the previous section.
To begin, we introduce/recall a bit of notation.
We denote by $q(x,\xi)$ and $\tilde{q}(x,\xi)$ the right and left eigenfunctions of the operator $L_{\xi}[u]$, respectively,
associated with the eigenvalue $\mu(\xi)$ bifurcating from the $(\xi,\mu(\xi))=(0,0)$ state.  Moreover, we assume
the normalization condition $\left<\tilde{q}(\cdot,\xi),q(\cdot,\xi)\right>_{L^2([0,1])}=1$.
Then we have the following representation for the Green kernel $G^I$.

\begin{lem}
Under assumptions (H1)-(H2) and (D1)-(D2), we have
\begin{align}
\left[G^I_{\xi}(x,t;y)\right]&=\phi(\xi)e^{\mu(\xi)t}q(x,\xi)\tilde{q}(y,\xi),\nonumber\\
G^I(x,t;y)&=\frac{1}{2\pi}\int_{\RM}e^{i\xi(x-y)}\left[G^I_{\xi}(x,t;y)\right]d\xi\label{GIrep}
\\ &
=\frac{1}{2\pi}\int_{\RM}e^{i\xi(x-y)}\phi(\xi)e^{\mu(\xi)t}q(x,\xi)\tilde{q}(y,\xi)d\xi\nonumber.
\end{align}
\end{lem}

\begin{proof}
The first equality is immediate from the spectral decomposition of elliptic operators on compact spatial domains.
Moreover, using the fact that the Fourier transform (either continuous or discrete) of the delta function is unity we have
\begin{align*}
\widehat{\delta_y}(\xi,x)&=\sum_{j\in\ZM}e^{2\pi ijx}\widehat{\delta_y}(\xi+2\pi j)
=\sum_{j\in\ZM}e^{2\pi ijx}e^{-i(\xi+2\pi j)y}\\
&=e^{-i\xi y}\sum_{j\in\ZM}e^{2\pi ij(x-y)}
=e^{-i\xi y}[\delta_y(x)].
\end{align*}
It follows that
\begin{align*}
G^I(x,t;y)&=\frac{1}{2\pi}\int_{-\pi}^\pi e^{i\xi x}\phi(\xi)P(\xi)e^{L_{\xi}[u]t}\widehat{\delta_y}(\xi,x)d\xi\\
&=\frac{1}{2\pi}\int_{-\pi}^\pi e^{i\xi (x-y)}\phi(\xi)P(\xi)e^{L_{\xi}[u]t}[\delta_y(x)]d\xi\\
&=\frac{1}{2\pi}\int_{-\pi}^\pi e^{i\xi(x-y)}[G^I_{\xi}(x,t;y)]d\xi,
\end{align*}
which yields the second equality by recalling that $\phi$ is supported in $[-\pi,\pi]$.
\end{proof}

In order to obtain nonlinear stability in the present context, it turns out that one can not simply use bounds as in the previous
section on the function $G^I$: see Remark \ref{rem:g1} below.  Instead, we must take extra care in obtaining our low-frequency linearized estimates by
separating out
%
the slow-decaying translation mode from the faster-decaying ``good" part of the solution operator.  To this
end, we define the function
\[
\tilde{e}(x,t;y):=\frac{1}{2\pi}\int_{\RM}e^{i\xi(x-y)}\phi(\xi)e^{\mu(\xi)t}\tilde{q}(y,\xi)d\xi
\]
and notice that
\[
G^I(x,t;y)=u'(x)\tilde{e}(x,t;y)+\frac{1}{2\pi}\int_{\RM}e^{i\xi(x-y)}\phi(\xi)e^{\mu(\xi)t}\left(q(x,\xi)-u'(x)\right)\tilde{q}(y,\xi)d\xi.
\]
By the analyticity of $q(x,\xi)$ in the variable $\xi$, we have that $q(x,\xi)-u'(x)=\mathcal{O}(|\xi|)$ for $|\xi|\ll 1$ and hence we expect
the difference $G^I(x,t;y)-u'(x)e(x,t;y)$ to decay {faster} than
the full Green kernel $G^I$.  This is the content of the
following proposition.

\begin{prop}\label{prop:refinedlfbds}
Under the hypothesis (H1)-(H2) and (D1)-(D2), the low-frequency Green kernel $G^I$ can be decomposed as
\[
G^I(x,t;y)=u'(x)\tilde{e}(x,t;y)+\widetilde{G}^I(x,t;y)
\]
where, for all $t>0$, and $2\leq p\leq\infty$ 
we have the estimate
\be\label{A}
\sup_y\left\|\widetilde{G}^I(\cdot,t;y)\right\|_{L^p(\RM)}
\lesssim(1+t)^{-\frac{1}{2}(1-1/p)-\frac{1}{2}}.
\ee
Moreover, for all $t>0$, $0\leq j,k,r, j+r\leq K+1$  we have the bound
\begin{equation}\label{eqn:ebds}
\sup_y\left\|\partial_x^j\partial_t^k\partial_y^r\tilde{e}(x,t;y)\right\|_{L^p(x;\RM)}\lesssim (1+t)^{-\frac{1}{2}(1-1/p)-\frac{(j+k)}{2}}.
\end{equation}
\end{prop}

\begin{proof}
Using the analyticity of the function $q(x,\xi)$ on the variable $\xi$, we have that
\[
\widetilde{G}^I(x,t;y)=\frac{1}{2\pi}\int_{\RM}e^{i\xi(x-y)}\phi(\xi)e^{\mu(\xi)t}\mathcal{O}(|\xi|)\tilde{q}(y,\xi)d\xi
\]
and hence the triangle inequality and (D2) yields
\begin{align*}
\sup_y\left\|\widetilde{G}^I(\cdot,t;y)\right\|_{L^\infty(\RM)}&\lesssim\left\| |\xi|e^{-\theta|\xi|^2t}\phi(\xi)\right\|_{L^1(\xi;\RM)}\\
&\lesssim(1+t)^{-1}.
\end{align*}
Moreover, noting that \eqref{GIrep} may be viewed itself as a Bloch-Fourier decomposition with respect
to the variable $z:=x-y$, with $y$ appearing as a parameter, we may use \eqref{eqn:parseval} to estimate
\begin{align*}
\sup_y \left\|\widetilde{G}^I(\cdot,t;y)\right\|_{L^2(\RM)}&\lesssim\sup_y\left\|\phi(\xi)e^{\mu(\xi)t}(q(x,\xi)-u'(x))\tilde{q}(y,\xi)\right\|_{L^2(x;L^2(\xi))}\\
&\lesssim\sup_y\left\|\phi(\xi)e^{-\theta|\xi|^2t}|\xi|\right\|_{L^2(\xi;[-\pi,\pi])}\sup_y\|\tilde{q}(y,\cdot)\|_{L^\infty(\xi;[-\pi,\pi])}\\
&\lesssim(1+t)^{-3/4},
\end{align*}
where we have used in a crucial way the boundedness of $\tilde{q}$.  By $L^p$ interpolation then, we obtain the desired
$L^p$ bounds on $\widetilde{G}^I(x,t;y)$.  Similar calculations yields the corresponding bound
on $\tilde{e}(x,t;y)$ by noting that $y$-derivatives do not improve decay while $x$- and $t$-derivatives improve decay
by a factor of $t^{-1/2}$ as above.
\end{proof}

\begin{remark}\label{rem:g1}
It is important to note that the Green kernel $G^I$ does not decay fast enough in this one-dimensional setting to
close our nonlinear iteration argument presented in the next section.  Indeed, using calculations as above (see also \cite{OZ})
one can verify that
\[
\sup_y\left\|G^I(\cdot,t;y)\right\|_{L^p(\RM)}\lesssim(1+t)^{-\frac{1}{2}(1-1/p)}
\]
for all $t>0$ and $2\leq p\leq \infty$.
Thus, by factoring out the translation mode from $G^I$ we gain an extra
$t^{-1/2}$ decay, which will end up being sufficient to close our nonlinear iteration arguments.
\end{remark}

We now combine the various bounds derived above to obtain estimates on the full Green kernel $G(x,t;y)$.
First off, let $\chi(t)$ be a smooth cut off function defined for $t\geq 0$
such that $\chi(t)=0$ for $0\leq t\leq 1$ and $\chi(t)=1$ for $t\geq 2$ and define
\[
e(x,t;y):=\chi(t)\tilde{e}(x,t;y).
\]

\begin{cor}\label{refinedgreenbds}
Under the hypothesis (H1)-(H2) and (D1)-(D2), the Green kernel $G(x,t;y)$ decomposes as
\[
G(x,t;y)=u'(x)e(x,t;y)+\widetilde{G}(x,t;y)
\]
where for all $t>0$, $2\leq p\leq \infty$, and
$0\leq j,k,r, j+r\leq K+1$,
 we have
\begin{align}\label{eqn:finalg}
\left\|\int_{-\infty}^\infty \widetilde{G}(\cdot,t;y)g(y)dy\right\|_{L^p(\RM)}&\lesssim(1+t)^{-\frac{3}{4}}
t^{-\frac{1}{2}(1/2-1/p)}\|g\|_{L^1(\RM)\cap L^2(\RM)}
\end{align}
and
\begin{align}\label{eqn:finale}
\left\|\int_{-\infty}^\infty\partial_x^j\partial_t^k\partial_y^r e(\cdot,t;y)g(y)dy\right\|_{L^p(\RM)}\lesssim(1+t)^{-\frac{1}{2}(1-1/p)-\frac{j+k}{2}}\|g\|_{L^1(\RM)}.
\end{align}
\end{cor}

\begin{proof}
The bound \eqref{eqn:finalg} follows immediately by considering the cases $0<t\leq 1$ and $t\geq 1$ separately.  Indeed,
using the $L^2\to L^p$ high-frequency bounds in Proposition \ref{prop:hfbds} for short time
and the $L^p$ bound \eqref{A} of Proposition \ref{prop:refinedlfbds}
together with the triangle inequality
 for large time yields the desired result.
The bound \eqref{eqn:finale} follows
similarly, using the $L^p$ bound
\eqref{eqn:ebds} of Proposition \ref{prop:refinedlfbds} together
with the triangle inequality.
\end{proof}

\section{Nonlinear Stability}\label{s:nstab}

With the above linearized estimates in hand, we are in a suitable position to prove nonlinear stability of the periodic traveling wave $\bar u(x)$
of the system of reaction diffusion equations \eqref{eqn:rd}.  Our main result is as follows.

\begin{theo}\label{main}
let $\bar u$ be a periodic standing-wave solution of \eqref{eqn:rd} and let $\tilde{u}(x,t)$
be any solution of \eqref{eqn:rd} such that $\|\tilde{u}-\bar u\|_{L^1(\RM)\cap H^K(\RM)}$ is sufficiently small.
Then assuming (H1)-(H2) and (D1)-(D2),
there exists a constant $C>0$ and a function $\psi(\cdot,t) \in W^{K,\infty}(\RM)$
such that for all $t\geq 0$ and $p\geq 2$ we have the estimates
\ba\label{eq:smallsest}
\|\tilde u-\bar u(\cdot -\psi(\cdot,t))\|_{L^p(\RM)}(t)&\le
C(1+t)^{-\frac{1}{2}(1-1/p)-\frac{1}{2}}
\|\tilde u(\cdot,t)-\bar u\|_{L^1\cap H^K(\RM)}|_{t=0},\\
\|\tilde u-\bar u(\cdot -\psi(\cdot,t))\|_{H^K(\RM)}(t)&\le
C(1+t)^{-\frac{3}{4}}
\|\tilde u(\cdot,t)-\bar u\|_{L^1(\RM)\cap H^K(\RM)}|_{t=0},\\
\|(\psi_t,\psi_x)(\cdot,t)\|_{H^{K}(\RM)}&\le
C(1+t)^{-\frac{3}{4}}
\|\tilde u(\cdot,t)-\bar u\|_{L^1(\RM)\cap H^K(\RM)}|_{t=0},\\
\ea
and
\begin{equation}\label{eq:stab1}
\begin{aligned}
\|\tilde u-\bar u\|_{ L^p(\RM)}(t), \; \|\psi(\cdot,t)\|_{L^p(\RM)}&\le
C(1+t)^{-\frac{1}{2}(1-1/p)} \|\tilde u(\cdot,t)-\bar u\|_{L^1(\RM)\cap H^K(\RM)}|_{t=0}.
\end{aligned}
\end{equation}
In particular, $\bar u$ is nonlinearly asymptotically $L^1\cap H^K\to H^K$ stable with estimate
\begin{equation}\label{eq:stab2}
\|\tilde u-\bar u\|_{ H^K(\RM)}(t), \; \|\psi(\cdot,t)\|_{H^K(\RM)}\le
C(1+t)^{-\frac{1}{4}} \|\tilde u(\cdot,t)-\bar u\|_{L^1(\RM)\cap H^K(\RM)}|_{t=0}
\end{equation}
for all $t\geq 0$.
\end{theo}

\br \label{compare}
\textup{
The regularity requirements stated here can be reduced to $f\in C^2$
and smallness of the initial data $v_0:=\tilde u-\bar u$ in $L^1\cap L^\infty$
as described in Remark \ref{regrmk}.
This is to be compared with Schneider's \cite{S} assumptions $f\in C^4$
and $v_0$ small in a weighted $H^{1/2+\delta}$ space ($\delta>0$)
bounding $\|(1+|x|^2)v_0\|_{L^\infty}$, hence $\|v_0\|_{L^1\cap L^\infty}$,
by Sobelev embedding.
}
\er

\subsection{Nonlinear Perturbation Equations}
Let $\tilde{u}(x,t)$ be a solution of the system of reaction diffusion equations
\[
u_t=u_{xx}+f(u)+cu_x
\]
and define $u(x,t)=\tilde{u}(x+\psi(x,t),t)$ for some unknown function $\psi:\RM^2\to\RM$
to be determined later.  Moreover, let $\bar{u}(x)$ be a stationary solution and define
\begin{equation}
v(x,t)=u(x,t)-\bar{u}(x)=\tilde{u}(x+\psi(x,t),t)-\bar{u}(x)\label{pertvar}
\end{equation}

\begin{lem}
For $v$, $u$ as above, we have
\ba\label{eqn:1nlper}
u_t-u_{xx}-f(u)-cu_x&=\left(\partial_t-L\right)\bar{u}'(x)\psi(x,t)
+\partial_x R\\
&+ (\partial_t-\partial_x^2)  S
+\left(f'(v(x,t)-\bar{u}(x))-df(\bar{u}(x))\right)\psi_x,
\ea
where
\[
R:= v\psi_t - v\psi_{xx}+  (\bar u_x(x) +v_x(x,t))\frac{\psi_x^2}{1+\psi_x}
= O\left(|v|(|\psi_t|+|\psi_{xx}|) +\left(\frac{|\bar u_x|+|v_x|}{1-|\psi_x|} \right)|\psi_x|^2\right)
\]
and
\[
S:=- v\psi_x =O\left(|v|\cdot|\psi_x|\right).
\]
\end{lem}

\begin{proof}
Using the fact that
$\tilde u_t-\tilde{u}_{xx}-f(\tilde{u})-c\tilde{u}_x=0$, it follows by a straightforward computation that
\begin{equation}\label{altform}
u_t+f(u)-u_{xx}-cu_{x}= \tilde u_x \psi_t-\tilde u_{t} \psi_x-(\tilde u_x \psi_x)_x+f(\tilde{u})\psi_x,
\end{equation}
where it is understood that the argument of the function $\tilde u$ and its derivatives appearing
on the righthand side are evaluated at $(x+\psi,t)$.  Moreover, by another direct calculation,
using the fact that
\[
L(\bar{u}'(x))=\left(\partial_x^2+c\partial_x+df(\bar{u})\right)\bar{u}'(x)=0
\]
by translation invariance, we have
\begin{equation*}
\left(\partial_t-L\right)\bar{u}'(x)\psi =\bar{u}_x\psi_t
-(\bar{u}_x\psi_{x})_{x}-(c\bar{u}_x+\bar{u}_{xx})\psi_x
=\bar{u}_x\psi_t -(\bar{u}_x\psi_{x})_{x}+df(\bar{u})\psi_x.
\end{equation*}
Subtracting, and using the facts that,
by differentiation of $(\bar u+ v)(x,t)= \tilde u(x+\psi,t)$,
\ba\label{keyderivs}
\bar u_x + v_x&= \tilde u_x(1+\psi_x),\\
\bar u_t + v_t&= \tilde u_t + \tilde u_x\psi_t,\\
\ea
so that
\ba\label{solvedderivs}
\tilde u_x-\bar u_x -v_x&=
-(\bar u_x+v_x) \frac{\psi_x}{1+\psi_x},\\
\tilde u_t-\bar u_t -v_t&=
-(\bar u_x+v_x) \frac{\psi_t}{1+\psi_x},\\
\ea
we obtain
\begin{align*}
u_t+ f(u) - u_{xx}&=
(\partial_t-L)\bar{u}'(x)\psi
+v_x\psi_t - v_t \psi_x - (v_x\psi_x)_x\\
&+ \left((\bar u_x +v_x)\frac{\psi_x^2}{1+\psi_x} \right)_x
 +\left(df(v+\bar{u})-df(\bar{u})\right)\psi_x,
\end{align*}
yielding \eqref{eqn:1nlper} by
$v_x\psi_t - v_t \psi_x = (v\psi_t)_x-(v\psi_x)_t$
and
$(v_x\psi_x)_x= (v\psi_x)_{xx} - (v\psi_{xx})_{x} $.
\end{proof}

\begin{cor}\label{cor:canest}
The nonlinear residual $v$ defined in \eqref{pertvar} satisfies
\be\label{veq}
\left(\partial_t-L\right)v=\left(\partial_t-L\right)\bar{u}'(x_1)\psi
+Q+ R_x +(\partial_x^2-\partial_t)S+T,
\ee
where
\be\label{eqn:Q}
Q:=f(v(x,t)+\bar{u}(x))-f(\bar{u}(x))-df(\bar{u}(x))v=\mathcal{O}(|v|^2),
\ee
\be\label{eqn:R}
R:= v\psi_t - v\psi_{xx}+  (\bar u_x +v_x)\frac{\psi_x^2}{1+\psi_x},
\ee
\be\label{eqn:S}
S:= v\psi_x =O(|v| |\psi_x|),
\ee
and
\be\label{eqn:T}
T:=\left(df(v+\bar{u})-df(\bar{u})\right)\psi_x=O(|v||\psi_x|).
\ee
\end{cor}

\begin{proof}
Straightforward Taylor expansion comparing \eqref{eqn:1nlper} and
$\bar u_t - f(\bar u)-\bar u_{xx}-cu_x=0$.
\end{proof}

\subsection{Integral Representation/$\psi$-Evolution Scheme}

Using Corollary \ref{cor:canest} and applying Duhamel's principle we obtain the integral (implicit) representation
\begin{align*}
v(x,t)=u'(x)&\psi(x,t)+\int_{-\infty}^\infty G(x,t;y)v_0(y)dy\\
&+\int_0^T\int_{-\infty}^\infty G(x,t-s;y)\left(Q+R_y+(\partial_y^2-\partial_s)S+T\right)(y,s)dy~ds
\end{align*}
for the nonlinear perturbation $v$.  Thus, if we define $\psi$ implicitly via the formula
\begin{align*}
\psi(x,t):=-&\int_{-\infty}^\infty e(x,t;y)v_0(y)dy\\
&-\int_0^t\int_{-\infty}^\infty e(x,t-s;y)\left(Q+R_y+(\partial_y^2-\partial_s)S+T\right)(y,s)dy~ds,
\end{align*}
we obtain the integral representation
\begin{equation}\label{eqn:vint}
\begin{aligned}
v(x,t)=&\int_{-\infty}^\infty \widetilde{G}(x,t;y)v_0(y)dy\\
&+\int_0^T\int_{-\infty}^\infty \widetilde{G}(x,t-s;y)\left(Q+R_y+(\partial_y^2-\partial_s)S+T\right)(y,s)dy~ds.
\end{aligned}
\end{equation}
Moreover, differentiating and recalling that $e(x,t;y)=0$ for $0<t\leq 1$ we obtain
\begin{equation}\label{eqn:psiint}
\begin{aligned}
\partial_t^k\partial_x^m\psi(x,t):=-&\int_{-\infty}^\infty \partial_t^k\partial_x^me(x,t;y)v_0(y)dy\\
&-\int_0^t\int_{-\infty}^\infty \partial_t^k\partial_x^m e(x,t-s;y)\left(Q+R_y+(\partial_y^2-\partial_s)S+T\right)(y,s)dy~ds.
\end{aligned}
\end{equation}
Together, these form a complete system in the variables $\left(v,\partial_t^k\psi,\partial_x^m\psi\right)$,
$0\leq k\leq 1$, $0\leq m\leq K+1$.  In particular, given a solution of the system we may afterward
recover the shift function $\psi$.

Now, from the original differential equation \eqref{veq} together with
\eqref{eqn:psiint}, we readily obtain
short-time existence and continuity with respect to $t$ of solution $(v,\psi_t,\psi_x)\in H^{K}$ by a standard
contraction-mapping argument treating the
linear $df(\bar u)v$ term of the lefthand side along with $Q,R,S,T,\psi \bar u'$
terms of the righthand side as sources in the heat equation.


\subsection{Nonlinear Iteration}

Associated with the solution $(u,\psi_t,\psi_x)$ of the integral system \eqref{eqn:vint}-\eqref{eqn:psiint}, we define
\begin{equation}\label{eqn:eta}
\eta(t):=\sup_{0\leq s\leq t}\|(v,\psi_t,\psi_x)\|_{H^{K}(x;\RM)}(s)(1+s)^{3/4}.
\end{equation}
By short time $H^K(\RM)$ existence theory, the quantity $\|(v,\psi_t,\psi_x)\|_{H^K(\RM)}$ is continuous so long
as it remains small.  Thus, $\eta$ is a continuous function of $t$ as long is it remains small.
We now use the linearized estimates of Section 3 to prove that if $\eta$ is initially small then it must remain so.

\begin{lem}\label{lem:eta}
For all $t\geq 0$ for which $\eta(t)$ is
sufficiently small, we have the estimate
\[
\eta(t)\leq C\left( E_0+\eta(t)^2\right)
\]
for some constant $C>0$, so long as $E_0:=\|v(\cdot,0)\|_{L^1(\RM)\cap H^{K}(\RM)}$ is also sufficiently small.
\end{lem}

\begin{proof}
To begin, notice that by the descriptions of $Q$, $T$, $R$, and $S$ in Corollary \ref{cor:canest} we have that
\begin{align*}
\|(Q,R_x,T)(\cdot,t)\|_{L^1(\RM)\cap L^2(\RM)}&\leq \|(v,v_x,\psi_t,\psi_x)\|_{L^1(x;\RM)}^2(t)+\|(v,v_x,\psi_t,\psi_x)\|_{L^2(x;\RM)}^2(t)\\
&\lesssim \eta(t)^2(1+t)^{-3/2}
\end{align*}
so long as $\|(v_x,\psi_x)(\cdot,t)\|_{L^\infty(\RM)}\leq \|(v,\psi_x)\|_{H^{K}(x;\RM)}(t)\leq \eta(t)$ remains bounded, and likewise
\[
\|(\partial_t-\partial_x^2)S(\cdot,t)\|_{L^1(\RM)\cap L^2(\RM)}\leq \|(v,\psi_x)\|_{H^1(x;\RM)}^2(t)+\|(v,\psi_x)\|_{H^{2}(x;\RM)}^2(t)\lesssim \eta(t)^2(1+t)^{-3/2}.
\]
Thus, applying the bound \eqref{eqn:finalg} of Corollary \ref{refinedgreenbds} to representations \eqref{eqn:vint}-\eqref{eqn:psiint}, we obtain for any
$2\leq p\leq\infty$ the bound
\begin{equation}\label{vbds}
\begin{aligned}
\|v(\cdot,t)\|_{L^p(\RM)}&\lesssim  (1+t)^{-\frac{1}{2}(1-1/p)-1/2}E_0\\
&~~~~~+\eta^2(t)\int_0^t(1+t-s)^{-3/4} (t-s)^{-\frac{1}{2}(1/2-1/p)}(1+s)^{-3/2}ds\\
&\lesssim\left(E_0+\eta(t)^2\right)(1+t)^{-\min\left(\frac{1}{2}(1-1/p)+\frac{1}{2},\frac{1}{2}(1-1/p)+1\right)}\\
&\lesssim\left(E_0+\eta(t)^2\right)(1+t)^{-\frac{1}{2}(1-1/p)-\frac{1}{2}}\\
\end{aligned}
\end{equation}
and similarly using \eqref{eqn:finale} we have
\begin{equation}\label{psibds}
\begin{aligned}
\|(\psi_t,\psi_x)(\cdot,t)\|_{W^{K+1,p}(\RM)}&\lesssim (1+t)^{-\frac{1}{2}(1-1/p)}E_0\\
&~~~~~+\eta(t)^2\int_0^t(1+t-s)^{-\frac{1}{2}(1-1/p)-\frac{1}{2}}(1+s)^{-3/2}ds\\
&\lesssim\left(E_0+\eta(t)^2\right)(1+t)^{-\frac{1}{2}(1-1/p)-\frac{1}{2}},
\end{aligned}
\end{equation}
yielding in particular that $\|(\psi_t,\psi_x)\|_{H^{K+1}}$
is arbitrarily small if $E_0$ and $\eta(t)$ are,%
\footnote{Note
that we have gained a necessary one degree of regularity in $\psi$, the regularity of $\psi$ being limited only by the regularity
of the coefficients of the underlying PDE \eqref{eqn:peq}.} thus verifying the hypothesis of Proposition \ref{prop:nd} below.
By the nonlinear damping estimate given
in Proposition \ref{prop:nd}, therefore,
 the size of $v$ in $H^K(\RM)$ can be controlled by
its size in $L^2(\RM)$ together with $H^K$ estimates on the
derivatives of the
phase function $\psi$.  In particular, we have for some positive constants $\theta_1$ and $\theta_2$
\begin{align*}
\|v(\cdot,t)\|_{H^K(\RM)}^2&\lesssim e^{-\theta_1 t}E_0^2+\left(E_0+\eta(t)^2\right)^2\int_0^t e^{-\theta_2(t-s)}(1+s)^{-3/2}ds\\
&\lesssim e^{-\theta_1 t}E_0^2+\left(E_0+\eta(t)^2\right)^2(1+t)^{-3/2}\\
&\lesssim \left(E_0+\eta(t)^2\right)^2(1+t)^{-3/2}.
\end{align*}
This estimate together with \eqref{psibds} in the case $p=2$ completes the proof.
\end{proof}

\begin{prop}\label{prop:nd}
Assuming $(H1)-(H2)$, let $v(\cdot,0)\in H^K(\RM)$ (for $v$ as in \eqref{pertvar}) and suppose that for some $T>0$ the $H^K(\RM)$
norm of $v$ and the $H^{K+1}(\RM)$ norms of $\psi_t(\cdot,t)$ and $\psi_x(\cdot,t)$ remain
bounded by a sufficiently small constant for all $0\leq t\leq T$.  Then there are constants
$\theta_1$ and $\theta_2$ such that
\[
\|v(\cdot,t)\|^2_{H^K(\RM)}\lesssim e^{-\theta_1 t}\|v(\cdot,0)\|_{H^K(\RM)}^2+\int_0^te^{-\theta_2(t-s)}\left(\|v(\cdot,s)\|^2_{L^2(\RM)}+\|(\psi_t, \psi_x)(\cdot,s)\|_{H^K(\RM)}^2\right)ds.
\]
for all $0\leq t\leq T$.  
\end{prop}

\begin{proof}
Subtracting from the equation \eqref{altform} for $u$
the equation for $\bar u$, we may write the
nonlinear perturbation equation as
\ba\label{vperturteq}
v_t - df(\bar u)v-v_{xx}-cv_x= Q
+ \tilde u_x \psi_t -\tilde u_{t} \psi_x - (\tilde u_x \psi_x)_x+f(\tilde{u})\psi_x,
\ea
where it is understood that derivatives of $\tilde u$ appearing
on the righthand side
are evaluated at $(x+\psi(x,t),t)$.
Using \eqref{solvedderivs} to replace $\tilde u_x$ and
$\tilde u_t$ respectively by
$\bar u_x + v_x -(\bar u_x+v_x) \frac{\psi_x}{1+\psi_x}$
and
$\bar u_t + v_t -(\bar u_x+v_x) \frac{\psi_t}{1+\psi_x}$,
and moving the resulting $v_t\psi_x$ term to the lefthand side
of \eqref{vperturteq}, we obtain
\ba\label{vperturteq2}
(1+\psi_x) v_t -v_{xx}-cv_x&=
-df(\bar u)v+ Q
+ (\bar u_x+v_x) \psi_t
\\ &\quad
- ((\bar u_x+v_x)  \psi_x)_x
+ \Big((\bar u_x+v_x) \frac{\psi_x^2}{1+\psi_x}\Big)_x+f(\tilde{u})\psi_x
\ea
Taking the $L^2$ inner product in $x$ of
$\sum_{j=0}^K \frac{(-1)^{j}\partial_x^{2j}v}{1+\psi_x}$
against (\ref{vperturteq2}), integrating by parts,
and rearranging the resulting terms,
we arrive at the inequality
\[
\partial_t \|v(\cdot,t)\|_{H^K(\RM)}^2 \leq -\theta \|\partial_x^{K+1} v(\cdot,t)\|_{L^2(\RM)}^2 +
C\left( \|v(\cdot,t)\|_{H^K(\RM)}^2
+
\|(\psi_t, \psi_x)(\cdot,s)\|_{H^K(\RM)}^2 \right)
,
\]
for some $\theta>0$, $C>0$, so long as $\|\tilde u\|_{H^K(\RM)}$ remains bounded,
and $\|v(\cdot,t)\|_{H^K(\RM)}$ and $\|(\psi_t, \psi_x)(\cdot,t)\|_{H^{K+1}(\RM)}$ remain sufficiently small.
Using the Sobolev interpolation
$
\|g\|_{H^K(\RM)}^2 \leq  \tilde{C}^{-1}\|\partial_x^{K+1} g\|_{L^2(\RM)}^2 + \tilde{C} \| g\|_{L^2(\RM)}^2
$
for $\tilde{C}>0$ sufficiently large, we obtain
\[
\partial_t \|v(\cdot,t)\|_{H^K(\RM)}^2(t) \leq -\tilde{\theta} \|v(\cdot,t)\|_{H^K(\RM)}^2 +
C\left( \|v(\cdot,t)\|_{L^2(\RM)}^2+\|(\psi_t, \psi_x)(\cdot,s)\|_{H^K(\RM)}^2 \right)
\]
from which the desired estimate follows by Gronwall's inequality.
\end{proof}


\noindent
\emph{Proof of Theorem \ref{main}:} Recalling that $\eta(t)$ is continuous so long as it remains small, it follows
by continuous induction using Lemma \ref{lem:eta}
 that $\eta(t)\le 2CE_0$ for all $t\ge 0$ provided that
$E_0 < 1/ 4C$ and (as holds without loss of generality) $C\ge 1$,
yielding by \eqref{eqn:eta} the result
\eqref{eq:smallsest} in the case $p=2$.  Similarly, using \eqref{vbds} and \eqref{psibds}, the result
\eqref{eq:smallsest} follows for any $2\leq p\leq \infty$ with uniform constant $C>0$.

Finally, notice that by \eqref{pertvar} we have
\[
\tilde{u}(x,t)-\bar{u}(x)=v(x-\psi(x,t),t)+\left(\bar{u}(x-\psi(x,t))-\bar{u}(x)\right)
\]
and hence the size of $\tilde{u}(x,t)-\bar{u}(x)$ in $L^p$ or $H^K$ is controlled by the corresponding
size of the function $(v+\psi)(x,t)$ in the respective norm.  Therefore, using \eqref{vbds} along with the estimate
\begin{align*}
\|\psi(\cdot,t)\|_{L^p(\RM)}&\lesssim E_0(1+t)^{-\frac{1}{2}(1-1/p)}+\eta(t)^2\int_0^t\left(1+t-s\right)^{-\frac{1}{2}(1-1/p)}(1+s)^{-\frac{3}{2}}ds\\
&\lesssim\left(E_0+\eta(t)^2\right)(1+t)^{-\frac{1}{2}(1-1/p)},
\end{align*}
which follows by
\eqref{eqn:finale}
for all $2\leq p\leq \infty$, we obtain \eqref{eq:stab1}.  Similarly
we obtain \eqref{eq:stab2}, which completes the proof.

\br\label{regrmk}
\textup{
Integrating by parts in representation formulae
\eqref{eqn:vint}--\eqref{eqn:psiint}
to exchange the $\partial_y$ and $(\partial_y^2-\partial_s)$
derivatives on $R$ and $S$ for $-\partial_y$ and
$(\partial_y^2+\partial_s)$ derivatives on
$\tilde G$, and noting that $\tilde G_y$ and
$ (\tilde G_{yy}+\tilde G_s)$ satisfy $L^p$ bounds equal
to
the bounds obtained on $ \tilde G$
times the time-integrable factor $(1+t)^{1/2}t^{-1/2}$,
we could alternatively close the iteration entirely within
$L^p$ spaces, avoiding the need for a nonlinear damping estimate
and obtaining $L^p$ stability, $p\ge 2$ with the same rates
given in Theorem \ref{main} for data merely small
in $L^1\cap L^\infty$.\footnote{
Here, we are using cancellation
in time and space derivatives, $ G_{yy}+G_s=df(\bar{u}) G + cG_y$.}
However, this is special to the semilinear case, whereas
the damping argument given here applies also to more general
quasilinear equations.
}
\er

\end{document}